\documentclass[%
english,%
twoside%,%
]{%
amsart%
}

\usepackage{amssymb}

\newcommand{\mfh}{{\mathfrak h}}
\newcommand{\mfH}{{\mathfrak H}}
\newcommand{\mfg}{{\mathfrak g}}

\newtheorem{corollary}{Corollary}
\newcommand{\bc}{\begin{corollary}}
\newcommand{\ec}{\end{corollary}}

\newtheorem{lemma}{Lemma}
\newcommand{\bl}{\begin{lemma}}
\newcommand{\el}{\end{lemma}}

\newcommand{\bp}{\begin{proof}} 
\newcommand{\ep}{\end{proof}}

\newtheorem{statement}{}
\newcommand{\bs}{\begin{statement}}
\newcommand{\es}{\end{statement}}

\newcommand{\unit}{\mathbf 1}
\newcommand{\C}{{\mathbb C}}

\newcommand{\R}{{\mathbb R}}
\newcommand{\Z}{{\mathbb Z}}
\newcommand{\opon}{\ltimes}

\DeclareMathOperator{\abs}{abs}

\DeclareMathOperator{\Aut}{Aut}

\DeclareMathOperator{\End}{End}

\DeclareMathOperator{\GL}{GL}

\DeclareMathOperator{\Hom}{Hom}

\DeclareMathOperator{\Id}{Id}

\DeclareMathOperator{\Ker}{Ker}

\DeclareMathOperator{\SU}{SU}

\title{%
	A cellular braid action and the Yang - Baxter equation
}
\date{%
	June 19, 1996
}
\author{%
	Mirko L\"{u}dde
} 
\address{%
	Institut f\"ur Reine Mathematik 
\\
	Humboldt Universit\"at zu Berlin 
\\
	Ziegelstrasse 13a 
\\
	D-10099 Berlin, Germany 
	}
\thanks{%
	Supported by {\sl Deutsche Forschungsgemeinschaft}%
}
\email{%
	luedde@mathematik.hu-berlin.de
}
\keywords{%
	braid group, Yang - Baxter equation, 
	Gauss - Manin and 
	Knizhnik - Zamolodchikov connections 
}
\begin{document} 
\begin{abstract} 
	Using a theorem of Schechtman - Varchenko 
	on integral expressions for 
	solutions of Knizhnik - Zamolodchikov equations 
	we prove that the solutions of the Yang - Baxter equation
	associated to complex simple Lie algebras 
	belong to the class of 
	generalised Magnus representations of the braid group. 
	Hence they can be obtained from the 
	homology of a certain cell complex, 
	or equivalently as group homology of 
	iterated free groups. 
\end{abstract} 
\maketitle  
\tableofcontents 

\bibliographystyle{%
amsalpha%
}

\section{%
	Introduction
}
The braid valued Burau representation 
$ B_n \hookrightarrow \GL _n(\Z [B_n \opon F_n]), $
where $ B_n $ is Artin's braid group on $ n $ strings 
and $ F_n $ is the free group of rank $ n, $
has been introduced
as a generalisation of the classical Burau representation 
in \cite{ConstantinescuLuedde1992}.
Here we give a construction of it 
in terms of elementary topology. 
Moreover we prove 
that the Yang - Baxter representations 
of the braid group 
associated to the 
quasitriangular Hopf algebras 
$ U_q(\mfg ), $ 
$ \mfg $ a complex simple Lie algebra, 
can be obtained from the braid valued Burau matrices
via a homological construction. 
Therefore these Yang - Baxter matrices 
belong to the class of generalised Magnus representations 
introduced in \cite{Luedde1992,Luedde19953}.

By a theorem of E.~Artin 
the braid group $ B_n $ 
faithfully acts on the free group $ F_n. $ 
This action can be understood 
by the lifting of paths in the 
Faddell - Neuwirth fibration of configuration spaces. 
In this way, 
an action of the braid group on certain cell complexes 
can be obtained as well.
This action is algebraically described by the 
braid valued Burau matrices.
The cellular action when projected to homology 
yields the monodromy of the 
natural flat connection (Gauss - Manin connection)
of homology fiber bundles associated 
to the Faddell - Neuwirth fibration. 
This might be compared with a theorem of Kohno 
on the reduced Burau representation 
describing the monodromy of the Jordan - Pochhammer 
local system \cite{Kohno19882}.

In \cite{SchechtmanVarchenko1991} 
there is constructed an imbedding of 
the flat Knizhnik - Zamolodchikov bundle 
associated to a Kac - Moody algebra
into the flat cohomology bundle. 
We combine this result 
with the Drinfeld - Kohno theorem, 
cf.~\cite{Kohno19882,Drinfeld1990,Kassel}, 
on the monodromy of the Knizhnik - Zamolodchikov connection.
This shows 
that the Yang - Baxter matrices associated to 
complex simple Lie algebras
can be obtained from the cellular braid action. 
Hence these Yang - Baxter braid representations 
can be understood in purely topological terms. 
To put it algebraically, 
they are quotients of tensor products 
of the braid valued Burau modules.

However, 
in order to obtain an explicit description 
of the Yang - Baxter matrices, 
one explicitly had to compute the homology 
(with local coefficients) 
of complexes, 
the investigation of which had been initiated by 
Vl.~Arnold, E.~Brieskorn, 
P.~Orlik, L.~Solomon and others, 
see \cite{OrlikTerao,SchechtmanVarchenko1991}.

Nevertheless the presented setting 
might serve as a step towards an understanding 
of the topological meaning of solutions of the 
Yang - Baxter equation and of link invariants 
constructed from them. 
A similar approach to the Jones polynomial 
was made in 
\cite{Lawrence1993,ConstantinescuLuedde1995}.
In particular, 
it would be interesting to know 
whether universal Vassiliev invariants, 
as constructed from a particular Knizhnik - Zamolodchikov 
connection, 
can be understood in this way, too.

The present work goes back to ideas in 
\cite{ChristeFlume1987,GomezSierra1991,RamirezRueggRuiz-Altaba1991} 
in the context of conformal field theory 
and it is a continuation of the investigation 
started in 
\cite{ConstantinescuLuedde1992,Luedde1992,Luedde19953}. 
Closely related work was done in \cite{Lawrence1990}.

{\bf Acknowledgement: }
I thank Professors F.~Constantinescu and R.~Flume 
for drawing my attention to \cite{SchechtmanVarchenko1991}.
I also thank 
Prof.~Th.~Friedrich and Klaus Mohnke 
for patient explanations, 
the participants of our seminar on knot theory 
for their interest 
and the 
'Graduiertenkolleg f\"ur Geometrie und Nichtlineare Analysis' 
at the Humboldt university at Berlin 
for financial support. 

\section{%
	The braid action on cell complexes
}
We compute the action of the braid group on a 
$ 1 $ - dimensional cell complex. 
From this action we obtain the braid valued Burau matrices, 
as defined (differently)
in \cite{ConstantinescuLuedde1992,Luedde1992,Luedde19953}. 
They describe the natural connection 
of homology bundles associated to 
the Faddell - Neuwirth fibration. 

Proofs of facts on the braid group that are used 
without further explanation can be found in 
\cite{Birman,BurdeZieschang}.

For every positive integer $ k $ 
let 
$ 
X_k := \{ (x_1, \ldots , x_k) \in (\R ^2) ^k; 
x_i \neq x_j \text{ if } i \neq j \} 
$ 
be the configuration space of $ k $ points 
in the euclidean plane. 
By a theorem of Faddell and Neuwirth, 
the projection 
$ p : X_{m+n} \rightarrow X_n $ 
is a topological fiber bundle 
with fibers 
$ 
Y_{x} := p^{-1}(x) = 
\{ (y_1, \ldots , y_m) \in X_m; y_i \neq x_j 
\text{ for } 
i \in \{ 1, \ldots , m \} , 
j \in \{ 1, \ldots , n \} 
\} . 
$ 
For $ m = 1 $ 
the fibers are $n$ - fold punctured euclidean planes. 

Given a point 
$ x \in X_n, $ 
we construct a subgroupoid 
of the fundamental groupoid 
of the base space $ X_n. $ 
The objects are the permuted points 
$ \{ \pi (x); \pi \in S_n \} $ 
($S_n$ the permutation group on $ n $ elements) 
and the morphisms are the path classes between these points. 
Since 
$ B_n \simeq \pi _1(X_n / S_n, x), $ 
there is a functor 
mapping this groupoid 
onto the braid group, 
by identifying all the objects of $ G. $  

By a general fact on fibrations 
there is another functor mapping the groupoid 
to the category of topological spaces and homotopy classes 
of maps. 
It sends a point 
$ x \in X_n $ 
to the fiber 
$ Y_{x} $ over the point 
and sends a path class 
in $ X_n $ from $ x $ to $ x' $ 
to a homotopy class of deformations 
of $ Y_{x} $ into $ Y_{x'}. $ 
These deformations can be chosen so as to 
map the base points of the fibers $ Y_{x} $ into each other. 

In the case $ m = 1 $ we have the following lemma. 
\bl[Braiding isotopies]
Let 
$
\gamma := (\gamma _1, \ldots , \gamma _n) : I \rightarrow X_{n}
$ 
($ I := [0,1] $) 
be a continuous path (a representative of a braid) 
between points $ x, x' $ with 
$
x = \pi (x'), 
$ 
$ \pi \in S_n $ 
a permutation.
\begin{enumerate}
\item
There is an isotopy 
$
\psi : I \times \R ^2 \rightarrow \R ^2 
$  
(a continuous family of homeomorphisms), 
connecting the identity 
$ \Id : \R ^2 \rightarrow \R ^2 $ 
with a homeomorphism 
$
\psi (1): \R ^{2} \rightarrow \R ^{2} 
$  
such that   
$
\gamma (t) = (\psi (t)(x_{1}), \ldots , \psi (t)(x_{n})). 
$ 
\item 
Any two such isotopies related to homotopic paths 
are homotopy equivalent.
\item
Let 
$ D \subset \R ^2 $ 
be a closed subset, homeomorphic to a disc 
and containing the set 
$ \{ \gamma _1(0), \ldots , \gamma _n(0) \} . $ 
Then the class of isotopies associated to $ \gamma $ 
contains a representative $ \psi $ such that 
$ 
\psi \mid (\R ^2 \backslash D) 
$
is the identity mapping.
\end{enumerate}
\el

We use this topological fact in order to construct 
a braid action on a cell complex 
and on its edge - path groupoid.

\bl[Braid action on a free groupoid in the plane]
Let $ G $ be a free groupoid defined by the set of objects 
$ 
\{ 0 + \epsilon, 1 \pm \epsilon, \ldots , n \pm \epsilon \}
$
for fixed $ 0 < \epsilon < 1 / 2, $ 
with base point 
$ P := 0 + \epsilon $ 
and with generating set of morphisms  
$ 
\{ w_0, \ldots , w_{n-1} \} 
\cup \{ l_1^{\pm }, \ldots , l_n^{\pm } \} , 
$
where 
$ w_k $ has domain 
$ k + \epsilon  $ 
and target 
$ k + 1 - \epsilon , $ 
such that 
$ 
w_k \in \Hom _G(k + \epsilon , k + 1 - \epsilon ), 
$ 
whereas 
$ l_i^{\pm} \in \Hom _G(i \mp \epsilon , i \pm \epsilon ). $ 
Then there is an imbedding
$ B_n \hookrightarrow \Aut (G,P) $ 
of the braid group $ B_n $ into 
the group of automorphisms 
of $ G $ preserving the base point.
Under this imbedding a braid generator $ \tau _i $ 
acts on the morphisms according to 
$$ 
\begin{array}{cccc}
\tau _i : 
& 
w_j \mapsto  
\left\{ \begin{array}{l}
w_j, 
\\
w_{i-1} l_i^+ w_i l_{1+i}^+,
\\
w_{i}^{-1},
\\
(l_{i}^-)^{-1} w_{i} (l_{1+i}^-)^{-1} w_{1+i},
\end{array} \right.
& 
l_j^{\pm } \mapsto  
\left\{ \begin{array}{l}
l_j^{\pm }, 
\\
l_j^{\pm }, 
\\
l_{1+i}^{\mp },
\\
l_{i}^{\mp }, 
\end{array} \right.
&
\begin{array}{l}
j \not \in \{ i-1, i, 1+i \} , 
\\
j = i - 1, 
\\
j = i, 
\\
j = 1 + i.
\end{array}
\end{array}
$$ 
\el
\bp 
We consider the $ 1 $ - dimensional cell complex imbedded 
into the complex plane. 
The $ 0 $ - cells are precisely the objects of $ G. $  
The $ 1 $ - cells are the positively oriented half circles 
$ 
l_k^{\pm } := k + S^{1}(\epsilon )^{\pm } = 
\{ z \in \C ; 
\abs (z - k) = \epsilon , \Im (z) \leq (\geq) 0 \} 
$ 
of radius $ \epsilon $ 
for $ k \in \{ 1, \ldots , n \} $ 
and the closed intervals 
$ 
w_m := \{ [m + \epsilon , 1 + m - \epsilon];  
m \in \{ 0, \ldots , n-1 \} \} 
$ 
of length $ 1 - 2 \epsilon . $ 
A braid generator $ \tau _i $ 
($ i \in \{ 1, \ldots , n-1 \} $) 
defines a class of isotopies 
of the plane. 
There is a representative isotopy mapping the points 
$ i $ and $ 1 + i $ into each other  
and being the identity outside an open disc around 
the point $ i + (1/2) $ of radius $ R $ with  
$ (1 / 2) + \epsilon < R < 1. $ 
The isotopy joins the identity homeomorphism 
with a self homeomorphism of the disc 
which can be chosen cellular. 
In this way the described action on the groupoid 
emerges. 
One can check that the braid relations 
$ \tau _i \tau _j = \tau _j \tau _i $ 
for 
$ \abs (i -j) \geq 2 $ 
and 
$ 
\tau _i \tau _{1+i} \tau _i = \tau _{1+i} \tau _i \tau _{1+i} 
$ 
hold. 
The map 
$ B_n \rightarrow \Aut (G,P) $ 
is a monomorphism, 
since the fundamental group 
$ \pi _1(G,P) \leq G $ 
is a free group $ F_n $ of rank $ n $ imbedded in $ G. $ 
Its generators are the loops 
$
f_j := \gamma _j l_j \gamma _j^{-1} 
$ 
around the points 
$ j \in \{ 1, \ldots , n \} $
with paths 
$ 
\gamma _j := w_0 (l_1^-)^{-1} w_1 (l_2^-)^{-1} \cdots w_{j-1}
$
from $ P $ to $ j - \epsilon $ 
and circles 
$ l_j := l_j^+ l_j^- $ 
around $ j. $
(There is a spanning tree consisting of all 
$ 1 $ - cells without the cells $ l_j^+. $ 
The paths $ \gamma _j $ are the unique paths in the 
spanning tree running from the base point 
to $ j - \epsilon. $ 
The cells $ l_j^+ $ 
therefore are in bijective correspondence 
with the generators of the fundamental group via 
$ 
l_j^+ \mapsto \gamma _j l_j^+ (\gamma _j (l_j^-)^{-1})^{-1}. 
$)
The action of the braid group restricted to this subgroup 
is precisely the action of Artin's theorem on the imbedding 
$ B_n \hookrightarrow \Aut (F_n). $
Indeed, we get 
$$
\begin{array}{cccc}
\tau _i : 
& 
\gamma _j \mapsto 
\left\{ \begin{array}{l}
	\gamma _j, 
\\
	\gamma _{j-1} l_{j-1}^+, 
\\
	\gamma _{j} l_{j}^+ w_j l_{1+j}^+, 
\end{array} \right. 
&
l_j \mapsto 
\left\{ \begin{array}{l}
	l_j, 
\\
	l_{j-1}^- l_{j-1}^+, 
\\
	l_{1+j}^- l_{1+j}^+, 
\end{array} \right.  
&
\begin{array}{l}
i \not \in \{ j - 1, j \} , 
\\
i = j - 1,
\\
i = j.
\end{array} 
\end{array}
$$ 
This combines to 
$$
\tau _i : f_j \mapsto 
\left\{ \begin{array}{ll}
	f_j, & i \not \in \{ j - 1, j \} , 
\\
	f_{j-1}, & i = j - 1,
\\
	f_j f_{1+j} f_j^{-1}, & i = j, 
\end{array} \right. 
$$ 
an action that is known to be faithful.
\ep 

Given a free connected groupoid $ (G,P) $ 
with base point $ P $ 
such that $ G $ is freely generated by a subset 
$  S \subset \cup _{x,y} \Hom _G(x,y) $ 
we can construct a free left $ \pi _1(G,P) $ module 
$ M. $ 
It is freely generated by elements $ [s] $ bijectively 
corresponding to the generators $ s \in S. $
For every object $ x $ of $ G $ we choose a path 
$ \gamma _x \in \Hom _G(P,x) $ 
joining the base point $ P $ with the point $ x. $
If 
$ E \leq \End (G,P) $ 
is a monoid of endomorphisms of $ G $ 
preserving the base point, 
the tensor product 
$ \Z [ E \opon \pi _1(G,P)] \otimes _{\pi _1(G,P)} M $ 
(the semidirect product 
$ E \opon \pi _1(G,P) $ 
is defined by the action 
of $ E $ on $ G $)
obtains a structure as right $ E $ module by setting 
$$ 
(1 \otimes [s]) e := 
e \otimes (e(\gamma _x) \cdot \gamma _{e(x)}^{-1}) \delta (e(s)), 
$$
where 
$ s \in \Hom _G(x,z), $ 
$ e(s) \in \Hom _G(e(x), e(z)), $
$ e \in E $ 
and $ \delta $ is the extended Fox derivation for 
groupoids. 
It is uniquely determined by the equations 
\begin{eqnarray*}
\delta (s) & = & [s], 
\\
\delta (g g') & = & 
\delta(g) + (\gamma _x \cdot g \cdot \gamma _y^{-1}) \delta (g'), 
\end{eqnarray*}
with 
$ g \in \Hom _G(x,y), $ 
$ g' \in \Hom _G(y,z). $ 
(Notice that this implies, 
$ \delta (\Id _x) =  0 $ 
and  
$ 
\delta (g^{-1}) = 
- \gamma _y \cdot g^- \cdot \gamma _x^{-1} \delta (g), $ 
analogous to the case of the common Fox derivation.)
The properties of $ \delta $ 
are precisely the properties of a lifting map,
which maps paths in the cell complex 
to the corresponding chains 
in the chain complex of the universal covering, 
see \cite{BurdeZieschang}.
Hence we have lifted the action of $ E $ 
on the groupoid $ G $ 
to an action on the chain complex of the universal covering 
of $ G. $

In this way, 
from the cellular action of the braid group 
on our cell complex
we obtain modules over the braid group.

\bl[Braid valued Burau modules from groupoid]
\begin{enumerate}
\item 
The 
$ \Z [B_n \opon F_n] $ - $ \Z [B_n] $ 
bimodule induced by the 
cellular braid action 
projects to a quotient 
carrying a braid representation in terms of the matrices 
(which map to the 
classical reduced Burau matrices in 
$ \GL _n(\Z [t, t^{-1}]) $ 
for 
$ \tau _i \mapsto 1, $ $ f_j \mapsto t $) 
$$
\tau _j
\left(
\begin{array}{ccccc}
\unit _{j-2}	&0	&0	&0	&0\\
0		&1	&f_j	&0	&0\\
0		&0	&-f_j	&0	&0\\
0		&0	&1	&1	&0\\
0		&0	&0	&0	&\unit _{n-j-1}
\end{array}
\right) 
$$ 
in 
$ \GL _n(\Z [B_n \opon F_n]) $
(where the semidirect product 
$ B_n \opon F_n $ 
is defined by the 
previously computed action of 
$ B_n $ onto $ \pi _1(G,P) \simeq F_n $  
and $ \unit _k $ is the $ k $ by $ k $ unit matrix).
\item 
It also contains a submodule 
with braid action described by the matrices
(which map to the classical unreduced Burau 
matrices) 
$$
\tau _i
\left(
\begin{array}{ccccc}
\unit _{i-1} & 0 & 0 & 0
\\
0 & 1 - f_{i} f_{i+1} f_{i}^{-1} & f_{i} & 0
\\
0 & 1 & 0 & 0
\\
0 & 0 & 0 & \unit _{n-i-1}
\end{array}
\right) .
$$
\end{enumerate}
\el
\bp
We choose paths 
$
\gamma _j^+ := 
w_0 (l_1^-)^{-1} w_1 (l_2^-)^{-1} \cdots w_{j-1} (l_{j}^-)^{-1}
$ 
from $ P $ to $ j + \epsilon $ 
and 
$
\gamma _j^- := 
w_0 (l_1^-)^{-1} w_1 (l_2^-)^{-1} \cdots w_{j-1} 
$ 
from $ P $ to $ j - \epsilon , $
$ j \in \{ 1, \ldots , n \} . $ 
These paths under the braid action transform as 
$$
\begin{array}{cccc}
\tau _i : 
& 
\gamma _j^- \mapsto 
\left\{ \begin{array}{l}
	\gamma _j^-, 
\\
	f_{j-1} \gamma _{j-1}^+,
\\
	f_j f_{1+j} \gamma _{1+j}^+,
\end{array} \right. 
&
\gamma _j^+ \mapsto 
\left\{ \begin{array}{l}
	\gamma _j^+,
\\
	\gamma _{j-1}^-,
\\
	f_j \gamma _{1+j}^-,
\end{array} \right.
&
\begin{array}{l}
i \not \in \{ j - 1, j \} , 
\\
i = j - 1,
\\
i = j.
\end{array}
\end{array}
$$ 
The representation that emerges on the free module 
is then given by 
\begin{eqnarray*}
\lefteqn{\tau _i : 1 \otimes _{\pi _1(G,P)} [w_j] \mapsto } & & 
\\
& & 
\tau _i \otimes _{\pi _1(G,P)} 
\left\{ \begin{array}{ll}
	\tau _i(\gamma _j^+) (\gamma _j^+)^{-1} \delta ( w_j ), 
		& j \not \in \{ i-1, i, 1+i \} , 
	\\
	\tau _i(\gamma _j^+) (\gamma _j^+)^{-1} 
	\delta ( w_{j} l_{1+j}^+ w_{1+j} l_{2+j}^+ ), 
		& j = i - 1, 
	\\
	\tau _i(\gamma _j^+) (\gamma _{1+j}^-)^{-1} \delta ( w_{j}^{-1} ), 
		& j = i, 
	\\
	\tau _i(\gamma _{j}^+) (\gamma _{j-1}^-)^{-1} 
	\delta ( (l_{j-1}^-)^{-1} w_{j-1} (l_{j}^-)^{-1} w_{j} ), 
		& j = 1 + i
\end{array} \right.
\end{eqnarray*}
which can be calculated to be  
\begin{eqnarray*}
\lefteqn{\tau _i : 1 \otimes _{\pi _1(G,P)} [w_j] \mapsto } & & 
\\
& & 
\tau _i \otimes _{\pi _1(G,P)} 
\left\{ \begin{array}{ll}
	[w_j], & j \not \in \{ i-1, i, 1+i \} , 
	\\
	{}[w_{j}] + [l_{1+j}^+] + f_{1+j} [w_{1+j}] + f_{1+j} [l_{2+j}^+],
		& j = i - 1, 
	\\
	- f_j [w_{j}], & j = i, 
	\\
	- [l_{j-1}^-] + [w_{j-1}] - [l_{j}^-] + [w_j], 
		& j = 1 + i.
\end{array} \right. 
\end{eqnarray*}
By a similar computation, 
we can find the transformation of the 
remaining generators.
It is clear already from our previous lemma, 
that the submodule of 
$ [l_j^{\pm }] $ cells is invariant under the braid action. 
We therefore can consider the quotient by this submodule 
and obtain a module of chains $ [w_j] $ alone, 
which yields the given matrix representation.
(This corresponds to the limit 
$ \epsilon \rightarrow 0, $ 
in which the length of the cells $ l_j^{\pm} $ 
approaches zero.)
The second matrix representation is obtained 
similarly by considering the braid action on the 
module generated by elements 
$ \delta(f_j) $ 
corresponding to the free group 
$ \pi _1(G,P). $ 
\ep

The braid valued Burau matrices can be iterated. 
There is an imbedding 
$ B_n \opon F_n \hookrightarrow B_{1+n} $ 
such that the matrices can be considered to 
have values in 
$ \GL _n(\Z B_{1+n}). $
Hence the matrix elements themselves can 
be represented by braid valued Burau matrices, 
and so on. 
The presented construction suggests 
that the $ m $ - fold iteration describes the braid 
action on $ m $ - dimensional cell complexes 
homotopic to the $ 2m $ - dimensional fibers 
$ Y_{x} $ of the fiber bundle 
$ p : X_{n+m} \rightarrow X_n. $ 
Using the fact that 
$ \pi _1(Y_{x}, y) \simeq F_n \opon \cdots \opon F_{n+m-1} $ 
(the semidirect products are defined by the action of 
the image of 
$ B_n \opon F_n \hookrightarrow B_{1+n} $ 
onto $ F_{1+n} $) 
and that all higher homotopy groups vanish, 
this indeed has been proven in \cite{Luedde19953}. 
More precisely, 
the (iterations of the) braid valued Burau matrices 
describe a braid action on a 
free resolution for 
$ \pi _1(Y_{x}, y) $ 
and therefore on the (group-) homology 
$ H_*(Y_x;L) \simeq  H_*(\pi _1(Y_x,y);L) $
for suitable local coefficient systems $ L. $

We summarise the content of the present section. 
\bs[Homology bundle with flat connection]
Let $ L $ be a locally constant sheaf 
of vector spaces on $ X_{n+m}. $
\begin{enumerate}
\item 
There is a vector bundle over 
$ X_n  $
with fibers 
$ H_*(Y_{x}, i_x^*(L)), $ 
possessing a natural flat connection
(Gauss - Manin connection in homology). 
\item 
The monodromy of the connection 
is described 
(on the level of cell complexes) 
by the (iterations of the) braid valued Burau modules. 
\end{enumerate}
\es
\bp
Consider the sheaf associated to the presheaf 
$ U \mapsto H_k(V, L \mid V), $ 
where $ U \subset X_n $ is any open set and 
$ V := p^{-1}(U) \subset X_{n+m}. $
For suitable neighborhoods $ U $ of $ x \in X_n, $ 
$ 
H_k (V, L \mid V) \simeq
H_k (U \times Y_{x}, L) \simeq
H_k (Y_{x}, i_x^*(L)), 
$ 
since $ L $ is locally constant. 
This determines the stalks of 
the sheaf to be 
$ H_k (Y_{x}, i_x^*(L)). $
It also shows, 
it is a locally constant sheaf of vector spaces. 
Hence it gives rise to a flat vector bundle. 
The second part of the assertion follows 
from our previous construction and remarks. 
\ep

The iteration procedure of 
the braid valued Burau matrices 
seems to produce 
a huge amount of braid representations 
(these representations yield the 
generalised Magnus representations of \cite{Luedde19953}), 
which, however, appear to be difficult to handle. 
Both facts will be illustrated in the next section, 
where two strong theorems are needed to show 
that we may obtain a lot of solutions of the Yang - Baxter equation 
from it. 

\section{%
	 Monodromy of Knizhnik - Zamolodchikov equations
}
In \cite{ChristeFlume1987,SchechtmanVarchenko1991} 
there have been given 
integral expressions for flat sections 
of the Knizhnik - Zamolodchikov vector bundles.
The families of homology cycles 
over which the integrals are taken 
are flat sections of the natural flat 
connection of the homology bundles 
of the Faddell - Neuwirth fibration.
Combined with the Drinfeld - Kohno theorem 
this shows that the Yang - Baxter matrices 
of complex simple Lie algebras can be 
understood topologically: 
they are induced from the cellular action of the braid 
group described by the braid valued Burau matrices. 

Let $ \mfg $ 
be a finite dimensional complex simple Lie algebra
with Cartan subalgebra $ \mfh \leq \mfg , $ 
with simple roots 
$ \alpha _i \in \mfh ^* $ 
and with Chevalley raising (lowering) 
operators $ e_i $ ($f_i$), $ i \in \{ 1 , \ldots , r \} , $ 
respectively. 
Let 
$ \Lambda _1, \ldots , \Lambda _n \in \mfh ^*, $ 
and let $ L(\Lambda _i) $ be 
the irreducible heighest weight module of weight $ \Lambda _i. $
Let 
$ C \in \mfg \otimes \mfg $ 
be the canonical element w.~r.~t.~the Killing form 
$ K : \mfg \otimes \mfg \rightarrow \C $ 
and let 
$ C_{i, j} $ 
be the representation of $ C $ on the $i$ - th and $j$ - th 
factor in 
$ L := L(\Lambda _1) \otimes \cdots \otimes L(\Lambda _n). $ 
Then the Knizhnik - Zamolodchikov equation 
with respect to 
$ \mfg , L $ 
and 
$ \kappa \in \C \backslash \{ 0 \} $ 
is the equation  
$
0 = 
d(s) - 
\kappa \sum _{1 \leq i < j \leq n} 
\frac{C_{i,j}(s) d(x_i - x_j)}{(x_i - x_j)}
$
for a local section 
$ s : U \rightarrow L $ 
of the trivial vector bundle 
$ L \times X_n \rightarrow X_n, $ 
see 
\cite{KnizhnikZamolodchikov1984,Kohno19882}. 

Let  
$ \lambda := \sum _{i=1}^r k_i \alpha _i \in \mfh ^* $
be a positive root, 
$ \sum k_i = m $ 
and let 
$ L_{\Lambda - \lambda} \subset L $
be the corresponding weight space, where 
$ \Lambda := \sum _{i=1}^{n} \Lambda _i. $
It has been shown in 
\cite{SchechtmanVarchenko1991} 
how to construct solutions of the equation with 
values in the subspace
$ 
V_{\lambda} := 
\cap _{i=1}^r \Ker 
(e_i : L_{\Lambda - \lambda} \rightarrow L).
$
For definiteness, we shortly review the construction. 

Consider the Faddell - Neuwirth fiber bundle 
$ p : Y \rightarrow X, $ 
$ Y := X_{m+n}, X := X_n, $ 
$ (y_1, \ldots , y_{m+n}) \mapsto (x_1, \ldots , x_n) := (y_1, \ldots , y_n) $
with fibers 
$ i_x : Y_{x} \hookrightarrow Y. $ 
We choose a basepoint 
$ y_0 := (1, \ldots , n+m) \in Y $ 
and identify 
$ \pi _1(Y,y_0) \simeq P_{n+m}, $ 
where $ P_{n+m} $ is the pure braid group on $ n+m $ 
strings. 
Its generators 
$ 
\{ \vartheta _{i,j} := \vartheta _{j,i}, 1 \leq i < j \leq n+m \} 
$ 
are represented by closed paths in $ Y, $ 
where the point $ j $ 
encircles the point $ i $ once in a positive orientation 
and does not encircle any other components of $ y_0. $

We define a homomorphism 
$ 
L \in \Hom (\pi _1(Y,y_0), \C \backslash \{ 0 \} ) 
$ 
by setting 
$ 
L(\vartheta _{i,j}) := \exp (2 {\sqrt -1} \pi a_{i,j}) 
$ 
for exponents 
$$ 
a_{i,j} := \kappa 
\left\{ \begin{array}{ll}
K(\Lambda _i, \Lambda _j), 
& 1 \leq i, j \leq n, 
\\
- K(\alpha _{\pi(i)}, \Lambda _j), 
& n+1 \leq i \leq n+m; 1 \leq j \leq n,  
\\
K(\alpha _{\pi(i)}, \alpha _{\pi(j)}), 
& n+1 \leq i,j \leq n+m.  
\end{array} \right. 
$$
Here we have set 
$ \pi : \{ 1, \ldots , m \} \rightarrow \{ 1 , \ldots , r \} $ 
such that 
$ \pi (j) = i $ 
if 
$ \sum _{p=1}^{i-1} k_p < j \leq \sum _{p=1}^{i} k_p. $ 
The local system equivalently is defined by the function 
$ 
l_{\lambda}
:= \prod _{1 \leq i < j \leq n+m} (y_i - y_j)^{a_{i,j}}.
$ 

Let 
$ (\Omega ^*_Y, d) $ 
and 
$ (\Omega ^*_X, d) $ 
be the complexes of sheaves of holomorphic 
differential forms on $ Y $ and $ X, $ 
respectively. 
Similarly, let 
$ \Omega ^*_L $ 
be the complex of sheaves of holomorphic 
differential forms on $ Y $ with values in 
the flat vector bundle (local system)
$ \C \times _L \tilde {Y} $
(where $ \tilde {Y} $ is the universal covering) 
determined by $ L. $
Furthermore, 
$ 
\Omega ^*_{L,p} := 
\Omega ^*_L / (p^*(\Omega ^1_X) \wedge \Omega ^*_L) 
$
is the quotient of the 
$ \Omega ^*_Y $ - module of $ L $ - valued forms 
by the submodule 
generated by the pulled back $ 1 $ - forms on $ X, $
i.e.~the relative de Rham complex 
of sheaves of $ p $ - vertical forms.
Finally, 
$ \Omega ^* := p_*(\Omega ^*_{L,p}) $ 
is the direct image sheaf on the base space $ X. $ 
Let 
$ \mfH ^* $ 
be the cohomology sheaf 
of the complex $ \Omega ^* $
and let 
$ \mfH _* $ 
be the sheaf of sections on $ X $ 
of the homology bundle induced by the 
Faddell - Neuwirth fibration 
and the dual local system $ L^*. $  
There is the de Rham map 
$ 
\int : \mfH _* \otimes (V_{\lambda} \otimes \mfH ^*)
\rightarrow V_{\lambda}, 
$ 
integrating vertical 
differential forms over homology classes of the fibers. 

In \cite{SchechtmanVarchenko1991} 
it was shown that there are 
holomorphic functions $ \phi _I $ on $ Y, $ 
being rational in $ \C ^{m+n} \supset Y, $
and vectors 
$ v_I \in L_{\Lambda - \lambda} $ 
(for $ I $ belonging to a suitable index set)
such that for 
$
\eta _I := (2 {\sqrt -1} \pi )^{-m}
\phi _I l_{\lambda}
dy_{1+n} \wedge \cdots \wedge  dy_{m+n} 
$
the function 
$ 
\eta := 
\sum _{I} 
\eta _I \otimes v_I 
$ 
has the following property
(which still holds for $ \mfg $ being a
Kac - Moody algebra of a symmetrisable generalised 
Cartan matrix).
\bs[Solutions of Knizhnik - Zamolodchikov equation]
Let 
$ \sigma \in \mfH _m (U) $ 
be a family of cycles that is 
flat with respect to the natural connection. 
Then the map 
$ U \rightarrow V_{\lambda}, $
$ 
x \mapsto \int _{\sigma (x)} \eta (x)
$ 
satisfies the Knizhnik - Zamolodchikov equation. 
For sufficiently small 
$ \abs (\kappa ) \neq 0 $ 
the set of these functions for varying $ \sigma $  
is a complete set of solutions. 
\es 
If $ \Lambda _1 = \ldots = \Lambda _n, $
the local system $ L $ defined above is invariant 
under the permutation of the coordinates 
$ (y_1, \ldots , y_n) $
of the base space. 
Therefore the holonomy of the natural connection in homology 
yields a representation of the braid group 
$ B_n \simeq \pi _1(X_n/S_n, x) $ 
on the cycles 
$ \sigma (x) \in H_m(Y_{x}; L^*). $ 
Knowing also the action of 
$ B_n $ 
on the function $ \eta , $ 
we obtain a braid representation 
from the monodromy of the integral expressions. 

\bc[K. - Z. monodromy from braid valued Burau matrices]
Let $ V $ be an irreducible module over 
the quasitriangular Hopf algebra $ U_q(\mfg ). $ 
Then the braid group representations 
$ B_n \rightarrow \Aut (V^{\otimes n}) $ 
induced from the universal R - matrix of $ U_q(\mfg ) $ 
can be obtained as group homology modules 
$ H_m(\pi _1(Y_{x}; y), i_x^*(L)) $ 
and with braid action induced from the braid valued 
Burau module. 
\ec
\bp
On the one hand, 
the monodromy of the Knizhnik - Zamolodchikov 
equation associated to $ V $ 
yields the representation of $ B_n $ 
coming from the universal R - matrix of $ U_q(\mfg ). $ 
This follows from the Drinfeld - Kohno theorem, 
see e.g~\cite[sect.~XIX.4, pp.~458]{Kassel}. 
On the other hand, 
the solutions of the Knizhnik - Zamolodchikov equation 
for small $ \kappa $ are given as integrals 
over cycles in  
$ H_m(Y_{x}; i_x^*(L)). $ 
This homology can be computed as group homology 
of $ \pi _1(Y_{x}, y_0) $
and it carries the braid representation 
induced from the braid valued Burau matrices 
\cite{Luedde19953}. 
Hence for small values of $ \kappa , $ 
both representations coincide. 
But from Chen's iterated integral expression 
for the monodromy of the Knizhnik - Zamolodchikov equation 
we know it is an entire function in the variable $ \kappa , $
see \cite{Kohno19882}. 
Similarly, 
the braid representation deduced from 
the braid valued Burau matrices 
is given by polynomial functions in 
$ \exp (2 {\sqrt -1} \pi \kappa ) , $ 
hence as well is entire in $ \kappa . $  
Therefore both representations 
coincide as functions of $ \kappa . $
\ep

\providecommand{\bysame}{\leavevmode\hbox to3em{\hrulefill}\thinspace}

\end{document}